\documentclass[12pt]{article}
\usepackage{amssymb,amsmath,mathrsfs}
\usepackage{hyperref}
\usepackage{graphicx}
\usepackage{color}
\usepackage[OT2,OT1]{fontenc}
\usepackage{upquote}
\usepackage{authblk}
\usepackage[numbers,sort&compress]{natbib}

\begin{document}

\title{\LARGE\bf A rational approximation of the Fourier transform by integration with exponential decay multiplier}

\bigskip
\author[1, 2]{\small Sanjar M. Abrarov}
\author[2, 3, 4]{\small Rehan Siddiqui}
\author[3, 4]{\small Rajinder K. Jagpal}
\author[1, 2, 4]{\small \\ Brendan M. Quine}

\affil[1]{\scriptsize Thoth Technology Inc., Algonquin Radio Observatory, Achray Road, RR6, Pembroke,~ON,~Canada,~K8A~6W7 \normalsize}
\affil[2]{\scriptsize Dept. Earth and Space Science and Engineering, York University, 4700 Keele St., Toronto,~ON,~Canada,~M3J~1P3 \normalsize}
\affil[3]{\scriptsize Epic College of Technology, 5670 McAdam Rd., Mississauga, ON, Canada, L4Z 1T2 \normalsize}
\affil[4]{\scriptsize Dept. Physics and Astronomy, York University, 4700 Keele St., Toronto, ON, Canada, M3J 1P3 \normalsize}

\date{October 16, 2022}
\maketitle
\maketitle

\begin{abstract}
Recently we have reported a new method of rational approximation of the sinc function obtained by sampling and the Fourier transforms. However, this method requires a trigonometric multiplier that originates from shifting property of the Fourier transform. In this work we show how to represent the Fourier transform of a function $f(t)$ in form of a ratio of two polynomials without any trigonometric multiplier. A MATLAB code showing algorithmic implementation of the proposed method for rational approximation of the Fourier transform is presented.
\vspace{0.25cm}
\vspace{0.25cm}
\end{abstract}

\section{Introduction}

The forward and inverse Fourier transforms of two related functions $f\left( t \right)$ and $F\left( \nu  \right)$ can be defined in a symmetric form as \cite{Bracewell2000, Hansen2014}
\begin{equation}\label{eq_1}
\mathcal{F}\left\{ {f\left( t \right)} \right\}\left( \nu  \right) = F\left( \nu  \right) = \int\limits_{ - \infty }^\infty  {f\left( t \right){e^{ - 2\pi i\nu t}}dt}
\end{equation}
and
\begin{equation}\label{eq_2}
{\mathcal{F}^{ - 1}}\left\{ {F\left( \nu  \right)} \right\}\left( t \right) = f\left( t \right) = \int\limits_{ - \infty }^\infty  {F\left( \nu  \right){e^{2\pi i\nu t}}d\nu },
\end{equation}
where variables $t$ and $\nu $ are the corresponding Fourier-transformed arguments in $t$-space and $\nu $-space, respectively (time $t$ vs. frequency $\nu $, for example).

Fourier transform methods are widely used in many applications including signal processing \cite{Bracewell2000, Hansen2014}, spectroscopy \cite{Goydaragh2021,Wang2021} and computational finance \cite{Colldeforns-Papiol2018, Bankole2019, Huang2022}.

There are several efficient methods have been reported for rational approximations in literature. For example, the rational approximations may be built on the basis of the Newman nodes \cite{Zhang2018}, Chebyshev nodes \cite{Brutman1997}, logarithmic nodes \cite{Fang2021} and so on.

Recently we have reported a new method of rational approximation of the Fourier transform \eqref{eq_1} as given by \cite{Abrarov2020}
\begin{equation}\label{eq_3}
\mathcal{F}\left\{ {f\left( t \right)} \right\}\left( \nu  \right) \approx {e^{2\pi i\nu a}}\sum\limits_{m = 1}^{{2^{M - 1}}} {\frac{{{A_m}\left( {\sigma + 2\pi i\nu} \right) + {B_m}}}{{C_m^2 + {{\left( {\sigma + 2\pi i\nu} \right)}^2}}}},
\end{equation}
where $M$ is an integer determining number of summation terms ${2^{M - 1}}$, $a$ is a shift constant, $\sigma $ is a decay (damping) constant and
$$
{A_m} = \frac{1}{{{2^{M - 1}}}}\sum\limits_{n = 0}^N { {f\left( {nh - a} \right){e^{\sigma nh}}\cos \left( C_m nh \right)}},
$$
$$
{B_m} = \frac{1}{{{2^{M - 1}}}}\sum\limits_{n = 0}^N { {f\left( {nh - a} \right){e^{\sigma nh}}{C_m}\sin \left( {{C_m nh}} \right)}},
$$
$$
{C_m} = \frac{{\pi \left( {2m - 1} \right)}}{{{2^M}h}}
$$
are expansion coefficients.

It has been noticed that approximation \eqref{eq_3} is not purely rational and there was a question whether or not a rational function of the Fourier transform \eqref{eq_1} in explicit form without any trigonometric multiplier of kind
\begin{equation}\label{eq_4}
{e^{2\pi i\nu a}} = \cos \left( {2\pi \nu a} \right) + i\sin \left( {2\pi \nu a} \right),
\end{equation}
depending on argument $\nu $, can be obtained \cite{Siddiqui2019-2020}. Theoretical analysis shows that this trigonometric multiplier originating from shifting property of the Fourier transform can be indeed excluded. As a further development of our work \cite{Abrarov2020}, in this paper we derive a rational function of the Fourier transform \eqref{eq_1} that has no any trigonometric multiplier of kind \eqref{eq_4}. Therefore, it can be used as an alternative to the Pad\'e approximation. To the best of our knowledge, this method of rational approximation of the Fourier transform \eqref{eq_1} for a non-periodic function $f\left(t\right)$ has never been reported in scientific literature.

\section{Derivation}

\subsection{Preliminaries}

Assume that $\operatorname{Re}\left\{f\left( t \right)\right\}$ is even while $\operatorname{Im}\left\{f\left( t \right)\right\}$ is odd such that $f:\mathbb{R}\to\mathbb{C}$, but $\operatorname{Re}\left\{f\right\}:\mathbb{R}\to\mathbb{R}$ and $\operatorname{Im}\left\{f\right\}:\mathbb{R}\to\mathbb{R}$. Then it is not difficult to see that the Fourier transform \eqref{eq_1} of the function $f\left( t \right)$ can be expanded into two integral terms as follows
\footnotesize
\[
\mathcal{F}\left\{ {f\left( t \right)} \right\}\left( \nu  \right) = F\left( \nu  \right) = 2\int\limits_0^\infty  {\operatorname{Re}\left\{f\left( t \right)\right\}\cos \left( {2\pi \nu t} \right)dt}  + 2\int\limits_0^\infty  {\operatorname{Im}\left\{f\left( t \right)\right\}\sin \left( {2\pi \nu t} \right)dt}.
\]
\normalsize
Assume also that the function $f\left( t \right)$ behaves in such a way that for some positive numbers ${\tau _1}$ and ${\tau _2}$ the following integrals
$$
{\int\limits_{{\tau _1}}^\infty  {\operatorname{Re}\left\{f\left( t \right)\right\}\cos \left( {2\pi \nu t} \right)dt} } \approx 0
$$
and
$$
{\int\limits_{{\tau _2}}^\infty  {\operatorname{Im}\left\{f\left( t \right)\right\}\sin \left( {2\pi \nu t} \right)dt} } \approx 0
$$
are negligibly small and can be ignored in computation. Consequently, we can approximate the Fourier transform as given by
\footnotesize
\begin{equation}\label{eq_5}
\mathcal{F}\left\{ {f\left( t \right)} \right\}\left( \nu  \right) = F\left( \nu  \right) \approx 2\int\limits_0^{{\tau _1}} {\operatorname{Re}\left\{f\left( t \right)\right\}\cos \left( {2\pi \nu t} \right)dt}  + 2\int\limits_0^{{\tau _2}} {\operatorname{Im}\left\{f\left( t \right)\right\}\sin \left( {2\pi \nu t} \right)dt}.
\end{equation}
\normalsize

Further the values $2{\tau _1}$ and $2{\tau _2}$ will be regarded as widths (pulse widths) for the real and imaginary parts of the function $f\left(t\right)$, respectively.

\subsection{New sampling method}

Consider a sampling formula (see, for example, equation (3) in \cite{Rybicki1989})
\begin{equation}\label{eq_6}
f\left( t \right) = \sum\limits_{n =  - N}^N {f\left( {{t_n}} \right){\rm{sinc}}} \left( {\frac{\pi }{h}\left( {t - {t_n}} \right)} \right) + \varepsilon \left( t \right),
\end{equation}
where
$$
{\rm{sinc}}\left( t \right) = \left\{ 
\begin{aligned}
&\frac{{\sin t}}{t}, \quad\,\, t \ne 0\\
&1, \qquad\quad t = 0,
\end{aligned} \right.
$$
is the sinc function, ${t_n}$ is a set of sampling points, $h$ is small adjustable parameter (step) and $\varepsilon \left( t \right)$ is error term. Fran\c{c}ois Vi\`{e}te discovered that the sinc function can be represented by cosine product\footnote{This equation is also attributed to Euler.} \cite{Kac1959, Gearhart1990}
\begin{equation}\label{eq_7}
{\rm{sinc}}\left( t \right) = \prod\limits_{m = 1}^\infty  {\cos \left( {\frac{t}{{{2^m}}}} \right)}.
\end{equation}

In our earlier publications we introduced a product-to-sum identity \cite{Quine2013}
\begin{equation}\label{eq_8}
\prod\limits_{m = 1}^M {\cos \left( {\frac{t}{{{2^m}}}} \right)}  = \frac{1}{{{2^{M - 1}}}}\sum\limits_{m = 1}^{{2^{M - 1}}} {\cos \left( {\frac{{2m - 1}}{{{2^M}}}t} \right)}
\end{equation}
and applied it for sampling \cite{Abrarov2015a, Abrarov2015b} as incomplete cosine expansion of the sinc function for efficient computation of the Voigt/complex error function. It is worth noting that this product-to-sum identity has also found some efficient applications in computational finance \cite{Ortiz-Gracia2016, Colldeforns-Papiol2018, Maree2018} involving numerical integration.

Comparing identities \eqref{eq_7} and \eqref{eq_8} immediately yields
$$
{\rm{sinc}}\left( t \right) = \mathop {\lim }\limits_{M \to \infty } \frac{1}{{{2^{M - 1}}}}\sum\limits_{m = 1}^{{2^{M - 1}}} {\cos \left( {\frac{{m - 1/2}}{{{2^{M - 1}}}}t} \right)} .
$$
Unlike equation \eqref{eq_7}, this limit consists of sum of cosines instead of product of cosines. As a result, its application provides significant flexibilities in various numerical integrations \cite{Abrarov2015a, Abrarov2015b, Ortiz-Gracia2016, Colldeforns-Papiol2018, Maree2018}.

Change of variable ${2^{M - 1}} \to M$ in the limit above leads to
$$
{\rm{sinc}}\left( t \right) = \mathop {\lim }\limits_{M \to \infty } \frac{1}{M}\sum\limits_{m = 1}^M {\cos \left( {\frac{{m - 1/2}}{M}t} \right)}.
$$
Therefore, by truncating integer $M$ and by making another change of variable $t \to \pi t/h$ we obtain
\begin{equation}\label{eq_9}
{\rm{sinc}}\left( {\frac{\pi }{h}t} \right) \approx \frac{1}{M}\sum\limits_{m = 1}^M {\cos \left( {\frac{{\pi \left( {m - 1/2} \right)}}{{M\,h}}t} \right)} , \qquad  - M\,h \leqslant t \leqslant M\,h.
\end{equation}
The right side of equation \eqref{eq_9} is periodic due to finite number of the summation terms. As a result, the approximation \eqref{eq_9} is valid only within the interval $t \in \left[ { - M\,h,M\,h} \right].$

At equidistantly separated sampling grid-points such that ${t_n} = nh$, the substitution of approximation \eqref{eq_9} into sampling formula \eqref{eq_6} gives
\footnotesize
\begin{equation}\label{eq_10}
f\left( t \right) \approx \frac{1}{M}\sum\limits_{m = 1}^M {\sum\limits_{n =  - N}^N {f\left( {nh} \right)\cos \left( {\frac{{\pi \left( {m - 1/2} \right)}}{{M\,h}}\left( {t - nh} \right)} \right)} } ,\qquad  - M\,h \leqslant t \leqslant M\,h.
\end{equation}
\normalsize
It is important that in sampling procedure the total number of the sampling grid-points $2N+1$ as well as the step $h$ should be properly chosen to insure that the widths $2\tau_1$ and $2\tau_2$ are entirely covered.

As we can see, the sampling formula \eqref{eq_10} is based on incomplete cosine expansion of the sinc function that was proposed in our previous works \cite{Abrarov2015a, Abrarov2015b} as a new approach for rapid and highly accurate computation of the Voigt/complex error function \cite{Armstrong1967, Abramowitz1972, Berk2013}. Computations we performed show that this method of sampling is particularly efficient in numerical integration.

\subsection{Even function}

Suppose that our objective is to approximate the sinc function ${\rm{sinc}}\left( {\pi \nu } \right)$. First we take the inverse Fourier transform \eqref{eq_2} of the sinc function
$$
{\mathcal{F}^{ - 1}}\left\{ {{\rm{sinc}}\left( \pi\nu  \right)} \right\}\left( t \right) = \int\limits_{ - \infty }^\infty  {{\rm{sinc}}\left( \pi\nu  \right){e^{ - 2\pi i\nu t}}d\nu }  = {\rm{rect}}\left( t \right),
$$
where
\[
{\rm{rect}}\left( t \right) = \left\{
\begin{aligned} & 1, \qquad {\rm{if}} \, \left|t\right| < 1/2
\\ & 1/2, \quad {\rm{if}} \, \left|t\right| =  1/2
\\ & 0, \qquad {\rm{if}} \, \left|t\right| > 1/2,
\end{aligned}  \right.
\]
is known as the rectangular function. This function is even since ${\rm{rect}}\left(t\right)={\rm{rect}}\left(-t\right)$. The rectangular function ${\rm{rect}}\left( t \right)$ has two discontinuities at $t =  - 1/2$ and $t = 1/2$. Therefore, it is somehow problematic to perform sampling over this function. However, we can use the fact that
\begin{equation}\label{eq_11}
{\rm{rect}}\left( t \right) = \mathop {\lim }\limits_{k \to \infty } \frac{1}{{\left( {2{t}} \right)^{2k} + 1}}.
\end{equation}
Thus, by taking a sufficiently large value for the integer $k$, say $k = 35$, we can approximate the rectangular function \eqref{eq_11} quite accurately as
\[
{\rm{rect}}\left( t \right) \approx f\left( t \right)=\frac{1}{{\left(2t\right)^{70} + 1}}.
\]

\begin{figure}[ht]
\begin{center}
\includegraphics[width=24pc]{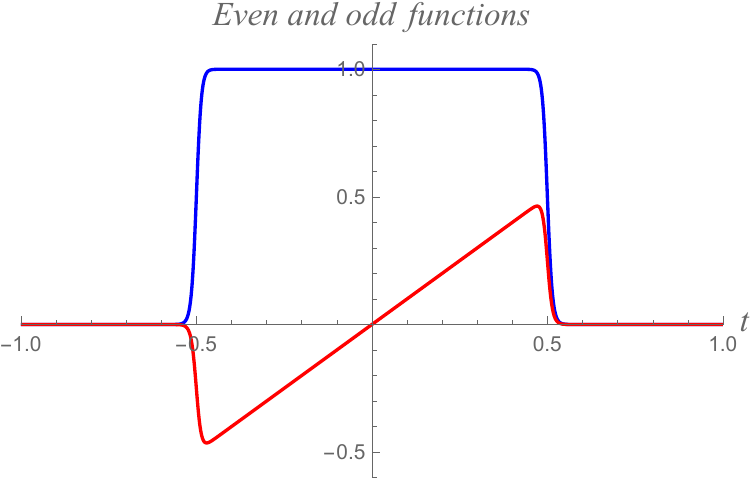}\hspace{2pc}%
\begin{minipage}[b]{28pc}
\vspace{0.3cm}
{\sffamily {\bf{Fig. 1.}} The even $1/\left(\left(2t \right)^{70}+1\right)$ and odd $t/\left(\left(2t \right)^{70}+1\right)$ functions shown by blue and red curves, respectively.}
\end{minipage}
\end{center}
\end{figure}

Figure 1 shows the function $f\left( t \right)=1/\left( {{{\left( {2t} \right)}^{70}} + 1} \right)$ by blue curve. As we can see from this figure, the function very rapidly decreases at $\left| t \right| > 1/2$ with increasing $t$. Therefore, we can take ${\tau _1} = 0.6$. Thus, the width of this function is $2{\tau _1} = 1.2$.

Sampling of function $f\left( t \right) = 1/\left( {{{\left( {2t} \right)}^{70}} + 1} \right)$ in accordance with equation \eqref{eq_10} results in a periodic dependence. Consequently, due to periodicity on the right side of equation \eqref{eq_10} it cannot be utilized for rational approximation of the Fourier transform. However, this problem can be effectively resolved by sampling the function $f\left( t \right){e^{\sigma t}}$ instead of $f\left( t \right)$ itself. This leads to
\footnotesize
\begin{equation}\label{eq_12}
\begin{aligned}
f\left( t \right){e^{\sigma t}} &\approx
\\ 
\frac{1}{M}&\sum\limits_{m = 1}^M {\sum\limits_{n =  - N}^N {f\left( {nh} \right){e^{\sigma nh}}\cos \left( {\frac{{\pi \left( {m - 1/2} \right)}}{{M\,h}}\left( {t - nh} \right)} \right)} } , \qquad  - M\,h \leqslant t \leqslant M\,h.
\end{aligned}
\end{equation}
\normalsize

\begin{figure}[ht]
\begin{center}
\includegraphics[width=24pc]{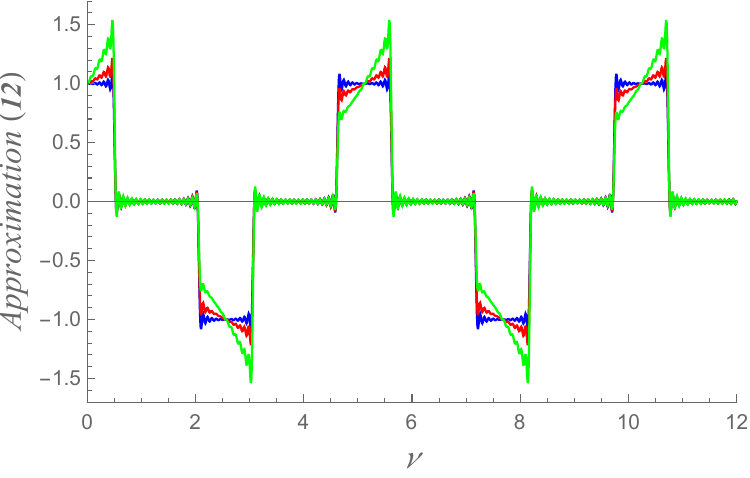}\hspace{2pc}%
\begin{minipage}[b]{28pc}
\vspace{0.3cm}
{\sffamily {\bf{Fig. 2.}}  Approximation \eqref{eq_12} to the function $f\left( t \right) e^{\sigma t} = e^{\sigma t}/\left[ {{{\left( 2t \right)}^{70}} + 1} \right]$ computed at $M = 32$, $N = 28$, $h=0.04$ with $\sigma = 0$ (blue curve), $\sigma = 0.25$ (red curve) and $\sigma = 0.75$ (green curve).}
\end{minipage}
\end{center}
\end{figure}

Figure 2 shows the results of computation for even function $f\left( t \right) = 1/\left( {{{\left( {2t} \right)}^{70}} + 1} \right)$ by approximation \eqref{eq_12} at $M = 32$, $N = 28$, $h = 0.04$ with $\sigma  = 0$ (blue curve), $\sigma  = 0.25$ (red curve) and $\sigma  = 0.75$ (green curve). As we can see from this figure, all three curves are periodic as expected. However, if the constant $\sigma $ is big enough, then slight rearrangement of equation \eqref{eq_12} in form
\begin{equation}\label{eq_13}
f\left( t \right) \approx \frac{{{e^{ - \sigma t}}}}{M}\sum\limits_{m = 1}^M {\sum\limits_{n =  - N}^N {f\left( {nh} \right){e^{\sigma nh}}\cos \left( {\frac{{\pi \left( {m - 1/2} \right)}}{{M\,h}}\left( {t - nh} \right)} \right)} } ,
\end{equation}
can effectively eliminate this periodicity due to presence of the exponential decay multiplier $e^{ - \sigma t}$ on the right side. This suppression effect can be seen from the Fig. 3 illustrating the results of computation for the even function $f\left( t \right) = 1/\left( {{{\left( {2t} \right)}^{70}} + 1} \right)$ by approximation \eqref{eq_13} at $M = 32$, $N = 28$  with $\sigma  = 0$ (blue curve), $\sigma  = 0.25$ (red curve) and $\sigma  = 0.75$ (green curve). As it is depicted by blue curve, at $\sigma  = 0$ the function is periodic. However, as decay coefficient $\sigma $ increases, the exponential multiplier ${e^{ - \sigma t}}$ suppresses all the peaks (except the first peak at the origin) such that the resultant function tends to become solitary along the entire positive $t$-axis. This tendency can be observed by red and green curves at $\sigma  = 0.25$ and $\sigma  = 0.75$, respectively. As a consequence, if the damping multiplier $\sigma $ is big enough, say greater than unity, the approximated function becomes practically solitary as the original function $f\left( t \right) = 1/\left( {{{\left( {2t} \right)}^{70}} + 1} \right)$ itself.

Thus, substituting approximation \eqref{eq_13} into equation \eqref{eq_5} and considering the fact that at sufficiently large $\sigma$ the function becomes solitary along positive $x$-axis, the upper limit $\tau_1$ of integration can be replaced by infinity as\footnote{For this integration we imply that the interval $2Nh$ along $t$-axis occupied by sampling grid-points is larger than the function width $2\tau_1 = 1.2$.}
\footnotesize
\[
\begin{aligned}
\mathcal{F}\left\{ {f\left( t \right)} \right\}&\left( \nu \right) =\mathcal{F}\left\{\operatorname{Re}\left\{{f\left( t \right)} \right\}\right\}\left( \nu \right) = \mathcal{F}\left\{ {\frac{1}{{{{\left( {2t} \right)}^{70}} + 1}}} \right\}\left( \nu  \right)
\\ 
\approx &\,2\int\limits_0^{\tau_{1}} {\left[ {\frac{{{e^{ - \sigma t}}}}{M}\sum\limits_{m = 1}^M {\sum\limits_{n =  - N}^N { \operatorname{Re}\left\{f\left( {nh} \right)\right\} {e^{\sigma nh}}\cos \left( {\frac{{\pi \left( {m - 1/2} \right)}}{{Mh}}\left( {t - nh} \right)} \right)} } } \right]{\cos\left(2\pi \nu t\right)}dt}
\\
\approx &\,2\int\limits_0^{\infty} {\left[ {\frac{{{e^{ - \sigma t}}}}{M}\sum\limits_{m = 1}^M {\sum\limits_{n =  - N}^N { \operatorname{Re}\left\{f\left( {nh} \right)\right\} {e^{\sigma nh}}\cos \left( {\frac{{\pi \left( {m - 1/2} \right)}}{{Mh}}\left( {t - nh} \right)} \right)} } } \right]{\cos\left(2\pi \nu t\right)}dt}.
\end{aligned}
\]
\normalsize
This integral can be taken analytically in form of rational function now and after some trivial rearrangements that exclude double summation, it follows that
\begin{equation}\label{eq_14}
\mathcal{F}\left\{\operatorname{Re}\left\{ {f\left( t \right)} \right\}\right\}\left( \nu  \right)\approx \sum\limits_{m = 1}^M {\frac{{{\alpha _m} + {\beta _m}{\nu ^2}}}{{{\kappa _m} + {\lambda _m}{\nu ^2} + {\nu ^4}}}},
\end{equation}
\begin{figure}[ht]
\begin{center}
\includegraphics[width=26pc]{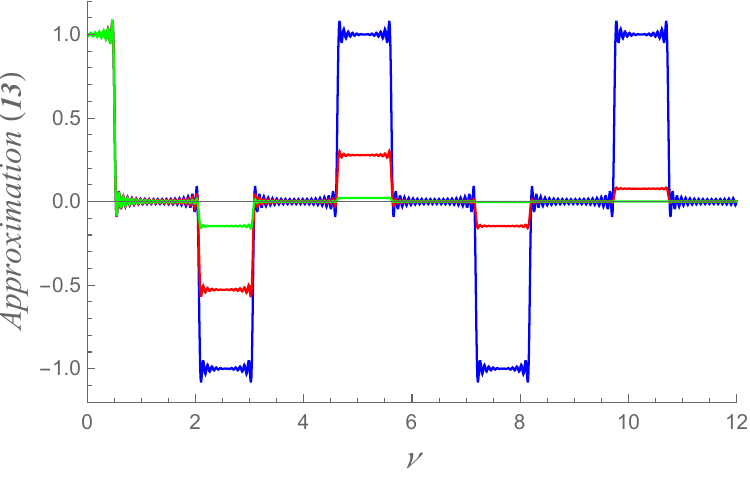}\hspace{2pc}%
\begin{minipage}[b]{28pc}
\vspace{0.3cm}
{\sffamily {\bf{Fig. 3.}} Evolution to the function $f\left( t \right) = 1/\left[ {{{\left(2t\right)}^{70}} + 1} \right]$ computed by approximation \eqref{eq_13} at $M = 32$, $N = 28$, $h=0.04$ with $\sigma = 0$ (blue curve), $\sigma = 0.25$ (red curve) and $\sigma = 0.75$ (green curve).}
\end{minipage}
\end{center}
\end{figure}
where the expansion coefficients are given by
\small
\[
{\alpha _m} = \frac{1}{8M\pi^4}\sum\limits_{n =  - N}^N \operatorname{Re}\left\{f\left( {nh} \right)\right\} e^{nh\sigma }\left( \mu_m^2 + \sigma ^2 \right)\left( \sigma \cos \left( nh \mu_m \right) + \mu_m\sin\left(nh\mu_m \right) \right),
\]
\normalsize
$$
\beta_m = \frac{1}{2 M \pi ^2}\sum\limits_{n =  - N}^N \operatorname{Re}\left\{f\left( {nh} \right)\right\} e^{nh\sigma} \left(\sigma\cos\left( nh\mu_m \right) - \mu_m\sin\left( nh \mu_m\right)\right),
$$
\[
\kappa_m = \frac{1}{16 \pi ^4}\left( \mu_m^2 + \sigma^2 \right)^2,
\]
\[
\lambda _m = \frac{1}{2\pi^2}\left( \sigma ^2 - \mu_m^2 \right)
\]
and
\[
\mu_m = \frac{\pi\left( m - 1/2 \right)}{Mh}.
\]

\begin{figure}[ht]
\begin{center}
\includegraphics[width=24pc]{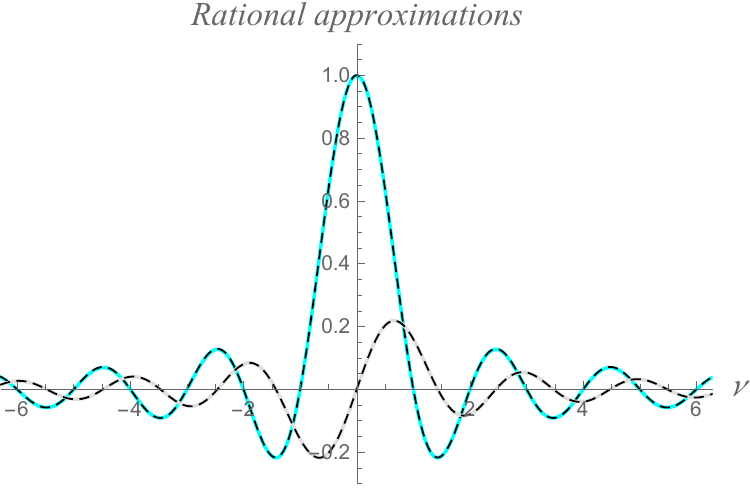}\hspace{2pc}%
\begin{minipage}[b]{28pc}
\vspace{0.3cm}
{\sffamily {\bf{Fig. 4.}} Approximations of the functions ${\rm{sinc}}\left(\pi\nu\right)$ and $\left({\sin \left( {\pi \nu } \right) - \pi \nu \cos \left( {\pi \nu } \right)}\right)/\left(2\left({\pi}{\nu}\right)^2\right)$  within interval $-2\pi\leq \nu \leq 2\pi$. Both approximations are obtained by equations \eqref{eq_14}, \eqref{eq_15} for input functions $1/\left(\left( 2t\right)^{70} +1\right)$ and $it/\left(\left( 2t\right)^{70} +1\right)$ at $M = 32$, $N = 28$, $h=0.04$, $\sigma = 2.7$ (light blue curve) and at $M = 32$, $N = 28$, $h=0.04$, $\sigma = 3$ (gray curve), respectively. The original functions ${\rm{sinc}}\left(\pi\nu\right)$ and $\left({\sin \left( {\pi \nu } \right) - \pi \nu \cos \left( {\pi \nu } \right)}\right)/\left(2\left({\pi}{\nu}\right)^2\right)$ are also shown by black dashed curves for comparison.}
\end{minipage}
\end{center}
\end{figure}

Figure 4 shows the original sinc function ${\rm{sinc}}\left( \nu  \right)$ and its approximation \eqref{eq_14} within the interval $ - 2\pi  \leqslant \nu  \leqslant 2\pi $ at $M = 32,$ $N = 28$, $h = 0.04$ and $\sigma  = 2.75$ by  black dashed and light blue curves, respectively. These two curves are not visually distinctive.

\subsection{Odd function}

Consider, as an example, the following function
$$
f\left( t \right) = \frac{it}{\left( 2t \right)^{70} + 1} \approx it\,{\rm{rect}}\left(t\right).
$$
We can see that the condition $t\,{\rm rect}\left( t \right)=-\left(-t\,{\rm rect}\left( -t \right)\right)$ for odd function in its imaginary part is satisfied. The function $\operatorname{Im}\left\{f\left( t \right)\right\}=t/\left( {{{\left( {2t} \right)}^{70}} + 1} \right)$ is shown in the Fig. 1 by red curve. We can take ${\tau _2} = 0.6$ and the width is $2{\tau _2} = 1.2$.

Using exactly same procedure as it has been described above and considering the fact that at sufficiently large $\sigma$ the upper limit $\tau_2$ of integration can be replaced by infinity, we can write\footnote{In this integration we imply again that the interval $2Nh$ along $t$-axis occupied by sampling grid-points is larger than the function width $2\tau_2 = 1.2$.}
\footnotesize
\[
\begin{aligned}
\mathcal{F}\left\{ {f\left( t \right)} \right\}&\left( \nu  \right) = \mathcal{F}\left\{i\,\operatorname{Im}\left\{{f\left( t \right)} \right\}\right\}\left( \nu  \right) =\mathcal{F}\left\{ i\,{\frac{t}{{{{\left( {2t} \right)}^{70}} + 1}}} \right\}\left( \nu  \right)
\\ 
\approx &\,2\int\limits_0^{\tau_{2}} {\left[ {\frac{{{e^{ - \sigma t}}}}{M}\sum\limits_{m = 1}^M {\sum\limits_{n =  - N}^N { \operatorname{Im}\left\{f\left( {nh} \right)\right\} {e^{\sigma nh}}\cos \left( {\frac{{\pi \left( {m - 1/2} \right)}}{{Mh}}\left( {t - nh} \right)} \right)} } } \right]{\sin\left(2\pi \nu t\right)}dt}
\\
\approx &\,2\int\limits_0^{\infty} {\left[ {\frac{{{e^{ - \sigma t}}}}{M}\sum\limits_{m = 1}^M {\sum\limits_{n =  - N}^N { \operatorname{Im}\left\{f\left( {nh} \right)\right\} {e^{\sigma nh}}\cos \left( {\frac{{\pi \left( {m - 1/2} \right)}}{{Mh}}\left( {t - nh} \right)} \right)} } } \right]{\sin\left(2\pi \nu t\right)}dt}.
\end{aligned}
\]
\normalsize
This leads to
\begin{equation}\label{eq_15}
\mathcal{F}\left\{i\,\operatorname{Im}\left\{{f\left( t \right)} \right\}\right\}\left( \nu  \right) \approx \sum\limits_{m = 1}^M {\frac{{{\eta_m}\nu  + {\theta _m}{\nu ^3}}}{{{\kappa _m} + {\lambda _m}{\nu ^2} + {\nu ^4}}}},
\end{equation}
where the expansion coefficients are
\small
\[
\eta_m = \frac{1}{4 M\pi ^3}\sum\limits_{n =  - N}^N \operatorname{Im}\left\{f\left( {nh} \right)\right\} e^{nh\sigma}\left( \left( \sigma ^2 - \mu_m^2 \right)\cos \left( nh \mu_m \right) + 2\sigma \mu_m\sin \left( nh\mu_m \right) \right)
\]
\normalsize
and
\[
\theta _m = \frac{1}{M\pi}\sum\limits_{n = - N}^N \operatorname{Im}\left\{f\left( {nh} \right)\right\} e^{nh\sigma}\cos\left( nh\mu_m \right).
\]

The Fourier transform of the function $it\,{\rm{rect}}\left(t\right)$ can be readily found analytically
$$
\begin{aligned}
\mathcal{F}\left\{ {it\,{\rm{rect}}\left( t \right)} \right\}\left( \nu  \right) =& \int\limits_{-\infty}^\infty  {it\,{\rm{rect}}\left( t \right){e^{ - 2\pi i\nu t}}dt}=\int\limits_{-1/2}^{1/2}  {it\,{\rm{rect}}\left( t \right){e^{ - 2\pi i\nu t}}dt}
\\
=&\int\limits_{-1/2}^{1/2}  {it\,{e^{ - 2\pi i\nu t}}dt} = \frac{{\sin \left( {\pi \nu } \right) - \pi \nu \cos \left( {\pi \nu } \right)}}{{2\left({\pi}{\nu}\right)^2}}.
\end{aligned}
$$

Gray curve in Fig. 4 illustrates the Fourier transform of the function $f\left( t \right)={it}/\left({\left(2t \right)^{70} + 1}\right)$ obtained by using approximation \eqref{eq_15} at $M=32$, $N=28$, $h=0.04$ and $\sigma=3$. The original function 
$$
\mathcal{F}\left\{ {it\,{\rm{rect}}\left( t \right)} \right\}\left( \nu  \right)=\frac{{\sin \left( {\pi \nu } \right) - \pi \nu \cos \left( {\pi \nu } \right)}}{{2\left({\pi}{\nu}\right)^2}}
$$
is also shown for comparison by black dashed curve. These two curves in the interval $ - 2\pi  \leqslant \nu  \leqslant 2\pi $ are also not distinctive visually.

\begin{figure}[ht]
\begin{center}
\includegraphics[width=24pc]{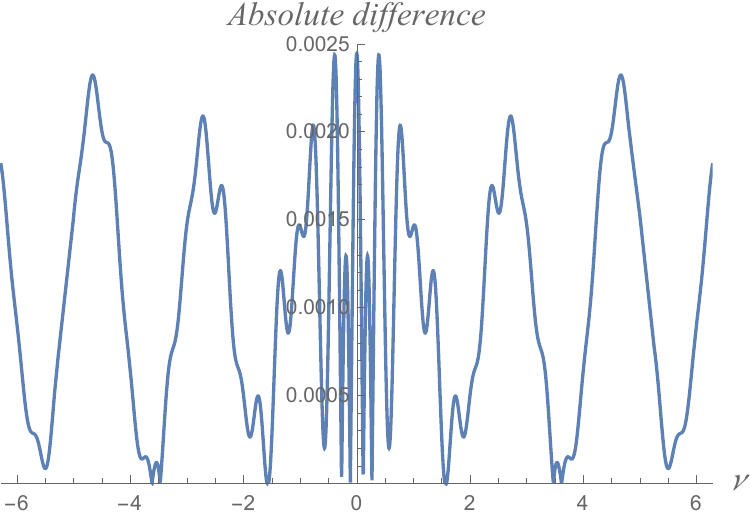}\hspace{2pc}%
\begin{minipage}[b]{28pc}
\vspace{0.3cm}
{\sffamily {\bf{Fig. 5.}} Absolute difference between the original sinc function $\rm{sinc} \left( \pi \nu \right)$ and its rational approximation \eqref{eq_14} for input function $f\left( t \right) = 1/\left(\left( 2t \right)^{70}+1\right)$ at $M = 32$, $N = 28$, $h=0.04$ and $\sigma = 2.7$.}
\end{minipage}
\end{center}
\end{figure}

\section{Accuracy}

Figures 5 shows the absolute difference between original sinc function ${\rm{sinc}}\left( {\pi \nu } \right)$ and its approximation \eqref{eq_14} for input function $f\left( t \right) = 1/\left( \left(2t\right)^{70}+1\right)$ calculated at $M=32$, $N=28$, $h=0.04$ and $\sigma = 2.7$. As we can see, the absolute difference within the interval $ - 2\pi  \leqslant \nu  \leqslant 2\pi $ does not exceed $2.5 \times 10^{-3}$. This accuracy is better than that of shown in our recent publication, where we used equation \eqref{eq_3} for the sinc function approximation (see Fig. 6 in \cite{Abrarov2020}).

\begin{figure}[ht]
\begin{center}
\includegraphics[width=24pc]{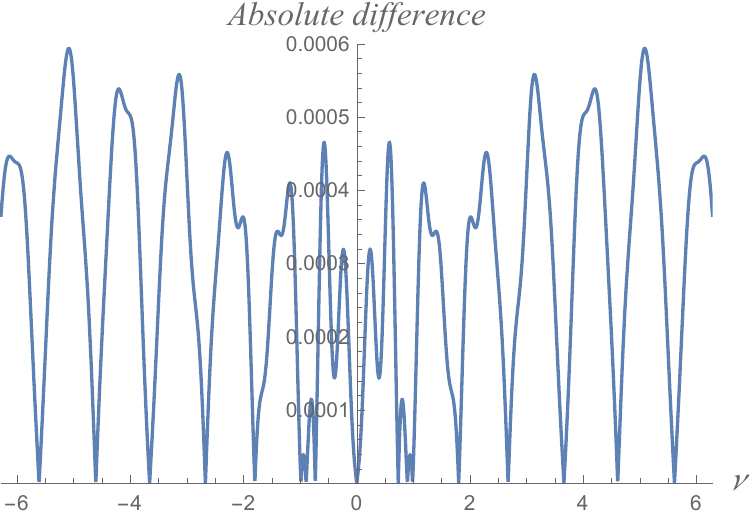}\hspace{2pc}%
\begin{minipage}[b]{28pc}
\vspace{0.3cm}
{\sffamily {\bf{Fig. 6.}} Absolute difference between the original function $\left({\sin \left( {\pi \nu } \right) - \pi \nu \cos \left( {\pi \nu } \right)}\right)/\left(2\left({\pi}{\nu}\right)^2\right)$ and its rational approximation \eqref{eq_15} for input function $f\left(t\right)=it/\left(\left( 2t \right)^{70} + 1 \right)$ at $M = 32$, $N = 28$, $h=0.04$ and $\sigma = 3$.}
\end{minipage}
\end{center}
\end{figure}

Figure 6 shows the absolute difference between original function given by $\left({\sin \left( {\pi \nu } \right) - \pi \nu \cos \left( {\pi \nu } \right)}\right)/\left(2\left({\pi}{\nu}\right)^2\right)$ and its approximation \eqref{eq_15} for input function $f\left( t \right) = it/\left( \left(2t\right)^{70}+1\right)$ calculated at $M=32$, $N=28$, $h=0.04$ and $\sigma = 3$. We can see that the absolute difference within the interval $ - 2\pi  \leqslant \nu  \leqslant 2\pi $ does not exceed $6\times 10^{-4}$.

\begin{figure}[ht]
\begin{center}
\includegraphics[width=24pc]{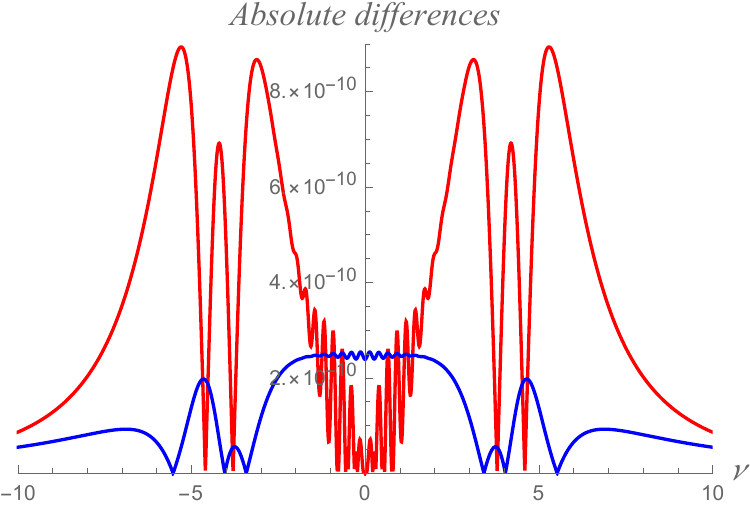}\hspace{2pc}%
\begin{minipage}[b]{28pc}
\vspace{0.3cm}
{\sffamily {\bf{Fig. 7.}} Absolute difference between the original functions $e^{-\nu^2}$ and its approximation \eqref{eq_14} for input function $\sqrt{\pi} e^{-\left(\pi t\right)^2}$ at $M = 16$, $N = 23$, $h=0.119$ and $\sigma = 6.9$ (blue curve). Absolute difference between the original functions $\nu e^{-\nu^2}$ and its approximation \eqref{eq_15} for input function $i\pi^{3/2} t e^{-\left(\pi t\right)^2}$ at $M = 16$, $N = 23$, $h=0.119$ and $\sigma = 5.9$ (red curve).}
\end{minipage}
\end{center}
\end{figure}

It should be noted that with more well-behaved functions we can obtain considerably higher accuracies. For example, 
suppose that we need to obtain the Fourier transform of the function $f\left( t \right)=\sqrt{\pi } e^{-(\pi  t)^2}+i\left[ \pi ^{3/2} t e^{-(\pi  t)^2 }\right]$ by using approximations \eqref{eq_14} and \eqref{eq_15}. Analytically, its Fourier transform is
$$
\mathcal{F}\left\{ \sqrt{\pi } e^{-(\pi  t)^2}+i\left[ \pi ^{3/2} t e^{-(\pi  t)^2 }\right] \right\}\left( \nu  \right)= e^{-\nu^2}+\nu e^{-\nu^2},
$$
where $e^{-\nu^2}$ is the Fourier transform of $\sqrt{\pi } e^{-(\pi  t)^2}$ while $\nu e^{-\nu^2}$ is the Fourier transform of $i\pi ^{3/2} t e^{-(\pi  t)^2}$. 

Blue curve in Fig. 7 corresponds to the absolute difference between function $e^{-\nu^2}$ and its approximation \eqref{eq_14} for input function $\sqrt{\pi } e^{-(\pi  t)^2}$ at $M=16$, $N=23$, $h=0.119$ and $\sigma=6.9$. Red curve in Fig. 7 corresponds to the absolute difference between function $\nu e^{-\nu^2}$ and its approximation \eqref{eq_15} for input function $ i\pi ^{3/2} t e^{-(\pi  t)^2}$ at $M=16$, $N=23$, $h=0.119$ and $\sigma=5.9$. We can see that with only 16 summation terms the absolute differences do not exceed $3\times 10^{-10}$ and $9\times 10^{-10}$. These results demonstrate that the rational approximations \eqref{eq_14} and \eqref{eq_15} can be highly accurate in the Fourier transform of well-behaved functions.

In our recent work \cite{Abrarov2018} we applied alternative method of sampling by using incomplete cosine expansion of the Gaussian function of kind $h\,e^{-\left(t/c\right)^2}/\left(c\sqrt{\pi}\right)$, where $c$ and $h$ are the fitting parameters. We have shown that this method of sampling can also be used to obtain high-accuracy computation of the Voigt/complex error function. In our future work will apply this method of sampling as an alternative that may reduce the absolute difference for rational approximations of the piecewise functions with discontinuities.

\section{Alternative representation}

For a function $f\left( t \right)=\operatorname{Re}\left\{f\left( t \right)\right\}+i\,\operatorname{Im}\left\{f\left( t \right)\right\}$, where its real part $\operatorname{Re}\left\{f\left( t \right)\right\}$ is even and its imaginary part $\operatorname{Im}\left\{f\left( t \right)\right\}$ is odd, we can write
$$
\mathcal{F}\left\{ {f\left( t \right)} \right\}\left( \nu  \right) \approx \sum\limits_{m = 1}^M {\frac{{{\alpha _m} + {\beta _m}{\nu ^2}}}{{{\kappa _m} + {\lambda _m}{\nu ^2} +  {\nu ^4}}}}  +\sum\limits_{m = 1}^M {\frac{{{\eta _m}\nu  + {\theta _m}{\nu ^3}}}{{{\kappa _m} + {\lambda _m}{\nu ^2} + {\nu ^4}}}}
$$
or
$$
\mathcal{F}\left\{ {f\left( t \right)} \right\}\left( \nu  \right) \approx \sum\limits_{m = 1}^M {\frac{{{\alpha _m} +{\eta_m}\nu  + {\beta _m}{\nu ^2} +{\theta _m}{\nu ^3}}}{{{\kappa _m} + {\lambda _m}{\nu ^2} + {\nu ^4}}}}.
$$
Using a Computer Algebra System (CAS) supporting symbolic programming it is not difficult to find coefficients ${p_k}$ and ${q_k}$ to represent this approximation as
\begin{equation}\label{eq_16}
\mathcal{F}\left\{ {f\left( t \right)} \right\}\left( \nu  \right) \approx \frac{{P\left( \nu  \right)}}{{Q\left( \nu  \right)}},
\end{equation}
where
$$
P\left( \nu  \right) = {p_0} + {p_1}\nu  + {p_2}{\nu ^2} +  \cdots  + {p_{4M - 2}}{\nu ^{4M - 2}}+ {p_{4M - 1}}{\nu ^{4M - 1}}
$$
and
\[
Q\left( \nu  \right) = {q_0} + {q_1}{\nu ^2} + {q_2}{\nu ^4} +  \cdots  + {q_{2M-1}}{\nu ^{4M-2}}+ {q_{2M}}{\nu ^{4M}},
\]
are polynomials of the orders $4M-1$ and $4M$, respectively.

Pad\'e approximation is one of the efficient methods to represent a function in form of ratio of two polynomials. Our preliminary numerical results show that the proposed new method of rational approximation may significantly extend the range $\left[\nu_{min},\nu_{max}\right]$ in coverage \cite{Siddiqui2019-2020} than the conventional Pad\'e approximation.

\section{MATLAB code and description}

The MATLAB code shown below is written as a function file {\sffamily{raft.m}} that can be simply copied and pasted to create {\sffamily{m}}-file in the MATLAB environment. The name of this function file originates from the abbreviation RAFT that stands for Rational Approximation of the Fourier Transform. The command {\sffamily{raft(opt)}} performs sampling and then computation of the expansion coefficients $\alpha_m$, $\beta_m$, $\eta_m$, $\theta_m$, $\kappa_m$, $\lambda_m$, $\mu_m$. Once the coefficients are determined, the program executes the Fourier transform according to equations \eqref{eq_14} and \eqref{eq_15} for even and odd input functions, respectively. The results of computations are generated in two plots. The first plot shows the Fourier transform of input function while the second plot illustrates its absolute difference.

There are four option values for {\sffamily{opt}} argument. At {\sffamily{opt = 0}}, {\sffamily{opt = 1}}, {\sffamily{opt = 2}} and {\sffamily{opt = 3}} the corresponding input functions are $\rm{rect}\left(t\right)$, $it\,\rm{rect}\left(t\right)$, $\sqrt{\pi}e^{-(\pi t)^2}$ and $i\pi^{3/2}te^{-(\pi t)^2}$. The default value is {\sffamily{opt = 0}} signifying that for the commands without argument {\sffamily{raft}} and {\sffamily{raft()}}, the value zero for {\sffamily{opt}} is assigned.

The authors did not attempt to optimize the code but rather to write it in a simple way with required comment lines in order to make it clear and intuitive for reading. The program was built and tested on MATLAB 2014a. However, the code should run in any version of MATLAB since it utilizes the most common commands.

\pagestyle{empty}
\bigskip
\footnotesize
\begin{verbatim}

function raft(opt)

% This function file performs the Rational Approximation of the Fourier 
% Transform (RAFT) for some functions and generates figures showing their
% Fourier transform (FT) and absolute differencies.
%
% SYNOPSIS:
%           opt = 0 provides FT for input function: f(t) = rect(t)
%           opt = 1 provides FT for input function: f(t) = i*t*rect(t)
%           opt = 2 provides FT for input function:
%                                   f(t) = sqrt(pi)*exp(-(pi*t)^2)
%           opt = 3 provides FT for input function:
%                                   f(t) = i*pi^(3/2)*t*exp(-(pi*t)^2)
% 
% The code is written by authors of this paper, York University, Toronto,
% Canada & Epic College of Technology, Mississauga, Canada.

clc
if nargin == 0
    disp('Missing input parameter! Option opt = 0 is assigned.')
    opt = 0;
end

switch opt
    case 0
        M = 32; N = 28; h = 0.04; sigma = 2.7; % define parameters
        n = -N:N; % define array n
        f = 1./((2*n*h).^70 + 1); % 1st even function
    case 1
        M = 32; N = 28; h = 0.04; sigma = 3; % define parameters
        n = -N:N; % define array n
        f = n*h./((2*n*h).^70 + 1); % odd function, imaginary unit i is ...
                                    % omitted
    case 2
        M = 16; N = 23; h = 0.119; sigma = 6.9; % define parameters
        n = -N:N; % define array n
        f = sqrt(pi)*exp(-(pi*n*h).^2); % 2nd even function
    case 3
        M = 16; N = 23; h = 0.119; sigma = 5.9; % define parameters
        n = -N:N; % define array n
        f = pi^(3/2)*n*h.*exp(-(pi*n*h).^2); % odd function, imaginary ...
                                             % unit i is omitted
    otherwise
        disp('Wrong parameter! Enter either 0, 1, 2 or 3.')
        return
end

% Initiate arrays for the expansion coefficients
alpha = zeros(M,1);
beta = zeros(M,1);
eta = zeros(M,1);
theta = zeros(M,1);
kappa = zeros(M,1);
lambda = zeros(M,1);

% Compute the expansion coefficients
f = f.*exp(n*h*sigma); % redefine the function array
for m = 1:M

    mu = pi*(m - 1/2)/M; % must be defined first

    if mod(opt,2) == 0 % if function is even, then:
        alpha(m) = 1/(8*M*pi^4)*sum(f.*((mu/h)^2 + sigma^2).*(sigma* ...
            cos(n*mu) + mu/(h)*sin(n*mu)));
        beta(m) = 1/(2*M*pi^2)*sum(f.*(sigma*cos(n*mu) - mu/h* ...
            sin(n*mu)));
    else % otherwise, if function is odd then:
        eta(m) = 1/(4*M*pi^3)*sum(f.*((sigma^2 - (mu/(h))^2)* ...
            cos(n*mu) + 2*sigma*mu/h*sin(n*mu)));
        theta(m) = 1/(M*pi)*sum(f.*cos(n*mu));
    end

    kappa(m) = 1/(16*pi^4)*((mu/h)^2 + sigma^2)^2;
    lambda(m) = 1/(2*pi^2)*(sigma^2 - (mu/h)^2);
end

% Rational approximation (14) in string format
strEven = ['f + (alpha(m) + beta(m)*nu.^2)./(kappa(m) + ', ...
    '    lambda(m)*nu.^2 + nu.^4)']; % for even function
% Rational approximation (15) in string format
strOdd = ['f + (eta(m)*nu + theta(m)*nu.^3)./(kappa(m) + ', ...
    '    lambda(m)*nu.^2 + nu.^4)']; % for odd function

nu = linspace(-2*pi,2*pi,1000); % define array for the argument nu        
f = 0; % reset function to zero
switch opt
    case {0, 2} % two even functions
        for m = 1:M % main computation for even functions
            f = eval(strEven); % use approximation (14)
        end

        if opt == 0 % assign some strings to plot 
            disp('Input function: f(t) = rect(t)')
            fStr = 'sinc(nu)';
            yStr = '\it{sinc(\pi\nu)}';
        else % if opt = 2, then:
            disp('Input function: f(t) = sqrt(pi)*exp(-(pi*t)^2)')
            fStr = 'exp(-nu.^2)';
            yStr = '\it{exp(-\nu^2)}';
        end
    case {1, 3} % two odd functions
        for m = 1:M % main computation for odd functions
            f = eval(strOdd); % use approximation (15)
        end

        if opt == 1 % assign some strings to plot
            disp('Input function: f(t) = i*t*rect(t)')
            fStr = '(sin(pi*nu) - pi*nu.*cos(pi*nu))./(2*(pi*nu).^2)';
            yStr = '\it{(sin(\pi\nu)-\pi\nu cos(\pi\nu))/(2(\pi\nu)^2)}';
        else % if opt = 3, then:
            disp('Input function: f(t) = i*pi^(3/2)*t*exp(-(pi*t)^2}')
            fStr = 'nu.*exp(-nu.^2)';
            yStr = '\it{\nu exp(-\nu^2)}';
        end
end

func2plot(nu,f,fStr,yStr); % call function to plot
absDiff(nu,f,fStr); % call to plot the abdolute difference

    function func2plot(nu,data2plot,funcStr,yAxisStr)

        % FIGURE 1
        figure1 = figure;
        axes1 = axes('Parent',figure1,'FontSize',12);
        xlim(axes1,[-2*pi,2*pi]);
        box(axes1,'on');
        grid(axes1,'on');
        hold(axes1,'all');
        plot1 = plot(nu,data2plot,nu,eval(funcStr),'Parent',axes1);
        set(plot1(1),'LineWidth',3,'Color',[0 1 1]);
        set(plot1(2),'LineStyle','--','LineWidth',2);
        xlabel('Parameter \it{\nu}','FontSize',14);
        ylabel(yAxisStr,'FontSize',14);
    end

    function absDiff(nu,data2plot,funcStr)

        % FIGURE 2
        figure2 = figure;
        axes2 = axes('Parent',figure2,'FontSize',12);
        xlim(axes2,[-2*pi,2*pi]);
        box(axes2,'on');
        grid(axes2,'on');
        hold(axes2,'all');
        plot2 = plot(nu,abs(eval(funcStr) - data2plot),'Parent',axes2);
        set(plot2(1),'LineWidth',1,'Color',[0 0 0]);
        xlabel('Parameter \it{\nu}','FontSize',14);
        ylabel('Absolute difference','FontSize',14);
    end
end

\end{verbatim}
\normalsize
\pagestyle{plain}

\section{Conclusion}

In this work we derived a rational approximation of the Fourier transform that with help of a CAS can be readily rearranged as
$$
\mathcal{F}\left\{ {f\left( t \right)} \right\}\left( \nu  \right) \approx \frac{{P\left( \nu  \right)}}{{Q\left( \nu  \right)}}.
$$
This method of the rational approximation is based on integration involving exponential decay multiplier ${e^{ - \sigma t}}$. The computational test we performed shows that this method of the Fourier transform can provide relatively accurate approximations of the Fourier transform even for the functions with discontinuities like ${\rm{rect}}\left( t \right)$ and $it\,{\rm{rect}}\left( t \right)$. Furthermore, this method shows that for the well-behaved function $f\left( t \right)=\sqrt{\pi } e^{-(\pi  t)^2}+i\left[ \pi ^{3/2} t e^{-(\pi  t)^2 }\right]$ with only $16$ summation terms the rational approximations \eqref{eq_14} and \eqref{eq_15} provide the Fourier transform with absolute differences not exceeding $3\times 10^{-10}$ and $9\times 10^{-10}$ for its real and imaginary parts, respectively. Our preliminary results indicate that the proposed method may be promising for rational approximation over the wide range $\left[\nu_{min},\nu_{max}\right]$.

\section*{Acknowledgments}

This work is supported by National Research Council Canada, Thoth Technology Inc., York University and Epic College of Technology.

\bigskip

\end{document}